\documentclass[12pt,a4]{amsart}
\usepackage{amsmath}
\usepackage{amssymb}
\usepackage{amsthm}
\usepackage{enumitem}
\usepackage{url}
\usepackage{soul}
\usepackage{times}
\usepackage[T1]{fontenc}

\usepackage{color}
\usepackage{dsfont}
\usepackage{epsfig}
\usepackage[hidelinks]{hyperref}
\usepackage{setspace}
\usepackage[authoryear,round,sort]{natbib}
\usepackage{comment}
\usepackage{booktabs}
\usepackage{float}
\usepackage{mathrsfs}
\usepackage{stmaryrd}
\numberwithin{equation}{section}

\usepackage{geometry}
\theoremstyle{plain}
\newtheorem{proposition}{Proposition}[section]

\newtheorem{theorem}{Theorem}[section]

\newtheorem{corollary}{Corollary}[section]

\theoremstyle{definition}

\theoremstyle{remark}
\newtheorem{rk}{Remark}[section]
\expandafter\let\expandafter\oldproof\csname\string\proof\endcsname
\let\oldendproof\endproof
\renewenvironment{proof}[1][\proofname]{%
  \oldproof[\noindent\textbf{#1.} ]%
}{\oldendproof}
\newcommand{\1}{\mathds{1}}

\newcommand{\be}{\begin{equation}}
\newcommand{\ee}{\end{equation}}
\newcommand{\by}{\begin{eqnarray*}}
\newcommand{\ey}{\end{eqnarray*}}

\renewcommand{\leq}{\leqslant}
\renewcommand{\geq}{\geqslant}
\usepackage{xcolor}
\definecolor{dark-red}{rgb}{0.4,0.15,0.15}
\definecolor{dark-blue}{rgb}{0.15,0.15,0.4}
\definecolor{medium-blue}{rgb}{0,0,0.5}
\hypersetup{
    colorlinks, linkcolor={red},
    citecolor={blue}, urlcolor={blue}
}
\allowdisplaybreaks

\begin{document}
\title{A Hoeffding's inequality for uniformly ergodic diffusion process}
\author{Michael C.H. Choi, Evelyn Li}\thanks{}
\address{Institute for Data and Decision Analytics, The Chinese University of Hong Kong, Shenzhen, Guangdong, 518172, P.R. China}
\email{michaelchoi@cuhk.edu.cn}
\address{School of Mathematics, Sun Yat-Sen University, Guangzhou, Guangdong, 510275, P.R. China}
\email{liylin6@mail2.sysu.edu.cn}
\date{\today}
\maketitle


\begin{abstract}

    In this note, we present a version of Hoeffding's inequality in a continuous-time setting, where the data stream comes from a uniformly ergodic diffusion process. Similar to the well-studied case of Hoeffding's inequality for discrete-time uniformly ergodic Markov chain, the proof relies on techniques ranging from martingale theory to classical Hoeffding's lemma as well as the notion of deviation kernel of diffusion process. We present two examples to illustrate our results. In the first example we consider large deviation probability on the occupation time of the Jacobi diffusion, a popular process used in modelling of exchange rates in mathematical finance, while in the second example we look at the exponential functional of a finite interval analogue of the Ornstein-Uhlenbeck process introduced by \cite{KS99}.
	\smallskip
	

\end{abstract}



\section{Introduction and main results}

The seminal work of \cite{H63}, which gives bound on large deviation probability of sum of bounded random variables, has now became one among many classical tools in probability theory. In particular, it has far reaching applications in statistics and machine learning, see for instance \cite{MR1} and the references therein. Hoeffding's inequality has then been refined and extended to various settings. For example, motivated by applications in Markov decision processes and reinforcement learning \cite{GO02} derives a Hoeffding's inequality for uniformly ergodic Markov chain, while \cite{B09} presents another method to prove Hoeffding's inequality in terms of the Drazin inverse of Markov chain. 

Inspired by the work cited above, we aim at extending Hoeffding's inequality to the setting of diffusion process. Contrary to the classical setting, we assume that the data stream arrives continuously from a uniformly ergodic diffusion process. The major difficulty of the analysis is then twofold. First, as we have a continuous data stream instead of discrete data points, previous analysis does not carry over to this setting easily. In addition, the dependency within the data stream complicates the situation. To overcome these difficulties, we employ classical martingale techniques for diffusion process as well as the notion of deviation kernel to aid our analysis. Comparing our result with the existing literature on concentration inequalities for diffusion processes \cite{GP07}, we argue that our proof is conceptually simpler since it utilizes similar techniques as in the discrete-time Markov chain case \cite{GO02,B09}. In addition, as we shall see in Corollary \ref{cor:main} below, it is readily applicable as long as we have the relevant eigenvalue information of the generator of the diffusion.

To this end, we fix our notation and introduce the tools we need for our main result Theorem \ref{thm:main} below. Let $(\Omega, \mathcal{F}, (\mathcal{F}_t)_{t \geq 0}, \mathbb{P})$ be a filtered probability space satisfying the usual conditions. Suppose that we have an ergodic diffusion process $X = (X_{t})_{t\geq 0}$ on state space $S$ with transition kernel $P$, transition density $P(x,dy)$ and stationary distribution $\pi$. We write $\mathbb{P}_{x}\left ( \cdot  \right ):= \mathbb{P}\left ( \cdot \mid X_{0}=x  \right )$ and $\mathbb{E}_{x}\left ( \cdot  \right ) :=
\mathbb{E}\left [ \cdot \mid X_{0}=x  \right ]$ to be the conditional probability and expectation when the process is initialized at $X_0 = x \in S$. $X$ is characterized by the infinitesimal generator $\mathcal{A}$, which acts on the space of twice differentiable functions and is defined to be
$$\mathcal{A} := \mu(x) \dfrac{d}{dx} + \dfrac{1}{2} \sigma^2(x) \dfrac{d^2}{dx^2},$$
where $\mu(x)$ and $\sigma^2(x)$ are respectively known as the drift and diffusion coefficient of $X$. A tool that we will use in the main result below is the deviation kernel $Q^{\sharp}$ of $X$, which is defined as 
$$Q^{\sharp}:=\int_{0}^{\infty } \left(P^t-\Pi\right) \,dt,$$
where $\Pi$ is the projection kernel with density $\Pi(x,dy) = \pi(dy)$ for all $x \in S$. It is well-known that the function $\hat{f} := Q^{\sharp} f$ solves the Poisson equation $-\mathcal{A}\hat{f}=f$, see e.g. \cite{GM96}. For bounded function $f$, we define the supremum norm to be $\left \| f \right \|:=\sup_x\left|f(x)\right|$. We also write $\left \| Q^{\sharp} \right\| := \sup_{\left \| f \right \| \leq 1} \left \| Q^{\sharp} f \right\|$ to be the induced operator norm of $Q^{\sharp}$ on the space of bounded functions. For further references on $Q^{\sharp}$, we refer readers to the work of \cite{CM15, Whitt92, Mao02}. On one-dimensional state space $S = (l,u)$, we now recall two fundamental notions associated with the diffusion $X$, namely the scale function $S(x)$ and the speed function $M(x)$. For $x \in \mathcal{X}$, these functions are defined by
\begin{align*}
S(x) := \int_{x_0}^x \exp \bigg\{ - \int_{x_0}^y \dfrac{2\mu(z)}{\sigma^2(z)}\,dz\bigg\}\, dy , \quad M(x) := \int_{l}^x \dfrac{2}{\sigma^2(y)} \exp \bigg\{ \int_{x_0}^y \dfrac{2\mu(z)}{\sigma^2(z)}\,dz\bigg\}\,dy,
\end{align*}
where $x_0 \in S$ is a fixed and arbitrary reference point. Their respective densities are given by
\begin{align*}
s(x) := \dfrac{d}{dx}S(x) = \exp \bigg\{ - \int_{x_0}^x \dfrac{2\mu(z)}{\sigma^2(z)}\,dz\bigg\}, \quad m(x) := \dfrac{d}{dx} M(x) = \dfrac{2}{\sigma^2(x) s(x)}.
\end{align*}

In this note, we are primarily interested in uniformly ergodic diffusions. That is, it is the class of ergodic diffusions such that the convergence to equilibrium in total variation distance is uniformly bounded by, for $t \geq 0$ and some constants $C < \infty$, $\beta > 0$,
$$\sup_{x \in S} ||P^t(x,\cdot) - \pi||_{TV} \leq C e^{-\beta t},$$
where $||P^t(x,\cdot) - \pi||_{TV} := \sup_{A} |P^t(x,A) - \pi(A)|$ is the total variation distance between $P^t(x,\cdot)$ and $\pi$. We write $\tau_y := \inf\{t \geq 0;~X_t = y\}$ to be the first hitting time of $y$ and $$t_{av}:=\int_{S\times S} \mathbb{E}_{x} \left[\tau _{y} \right]\pi(dx)\pi(dy)$$
to be the average hitting time of $X$. While verifying uniform ergodicity can be quite difficult, it turns out that, according to \cite[Theorem 2.2]{CM15}, uniform ergodicity for diffusion on $(0,u)$ with reflecting boundary at $0$ is equivalent to a few readily checkable conditions on $t_{av}$, $s(x)$ and $m(x)$:
\begin{proposition}[Necessary and sufficient conditions for uniform ergodicity \cite{CM15}]\label{prop:nscue}
	Given a ergodic diffusion $X$ on $S = (0,u)$ with reflecting boundary at $0$, the following statements are equivalent:
	\begin{enumerate}
		\item $X$ is uniformly ergodic;
		\item\label{it:S} $\int_S m([x,u]) s(x)\, dx < \infty$;
		\item $t_{av} < \infty$;
		\item\label{it:E} $\sigma_{ess}(\mathcal{A}) = \emptyset$ and $$\sum_{i \geq 1} \dfrac{1}{\lambda_i} < \infty,$$
		where $\sigma_{ess}(\mathcal{A})$ is the essential spectrum of $\mathcal{A}$ and $(\lambda_i)_{i \geq 1}$ are the non-zero eigenvalues of $-\mathcal{A}$.
	\end{enumerate}
\end{proposition}
At times it maybe easier to check item \eqref{it:S} as it depends on $\mu(x)$ and $\sigma^2(x)$ through $s(x)$ and $m(x)$, while at other times when eigenvalues information are available perhaps it is more convenient to check item \eqref{it:E}. As a simple illustration of item \eqref{it:S}, we consider the class of diffusions with $\mu(x) = 0, \sigma^2(x) = 2(1+x)^{\gamma}$ and $S = (0, \infty)$, where $\gamma > 2$ is a parameter. This class is first studied in \cite{Mao02}. It is easy to see that $s(x) = 1$ and $m(x) = (1+x)^{-\gamma}$. As a result, item \eqref{it:S} now reads
$$\int_S m([x,\infty]) s(x)\, dx = \dfrac{1}{(\gamma-1)(\gamma-2)} < \infty,$$ 
and so this class of diffusions with $\gamma >2$ are uniformly ergodic. For illustration of item \eqref{it:E}, we defer the readers to Corollary \ref{cor:main} when we discuss the Jacobi process. In view of Proposition \ref{prop:nscue}, for uniformly ergodic $X$ we have
$$||Q^{\sharp}|| \leq 2 t_{av} < \infty,$$
where the first inequality follows from \cite[Theorem $1.1$]{Choi18}. In other words, for uniformly ergodic diffusion the induced operator norm $||Q^{\sharp}||$ is finite. Such a term will appear in our version of Hoeffding's inequality Theorem \ref{thm:main} below.

With the above notation, we are now ready to state our main result. It follows from the classical ergodic theorem that for bounded function $f$, the time average $\frac{1}{t}\int_{0}^{t}f(X_{s})ds$ converges almost surely to the space average $\pi(f) := \int_{S}f(x) \pi(dx)$ as $t \to \infty$, see e.g. \cite[Theorem $12.2$]{BW09}. In our main result below, we present non-asymptotic probabilistic error bound of such convergence:

\begin{theorem}\label{thm:main}
	Suppose that $X = \left(X_{t} \right)_{t\geq 0}$ is a uniformly ergodic diffusion and $f$ is a bounded function. For any $\varepsilon > 0$, $t > \frac{2\left \| f \right \|\left \| Q^{\sharp} \right \|}{\varepsilon}$, $x \in S$,
	$$\mathbb{P}_{x}\left( \frac{1}{t}\int_{0}^{t}f(X_{s})ds-\pi(f)\geqslant \varepsilon  \right) \leqslant \exp\left \{ \frac{-2\left ( t\varepsilon -2\left \| f \right \|\left \|  Q^{\sharp} \right \| \right )^{2}}{(t+1)\left \| f \right \|^{2}\left ( 2\left \|  Q^{\sharp} \right \| +1\right )^{2}} \right \},$$
	where $Q^{\sharp}$ is the deviation kernel of the process $X$.
\end{theorem}

\begin{rk}[On the assumption of bounded $f$]
	As usual in the Hoeffding's inequality literature, our main result Theorem \ref{thm:main} requires the function $f$ to be bounded. This assumption is crucial when we apply the classical Hoeffding's lemma \cite[Lemma $8.1$]{MR1} to certain martingale difference sequence in \eqref{eq:part1} and \eqref{eq:part2} below, which only holds when the random variable of interest is bounded. Although there is extension of the  Hoeffding's lemma to non-negative random variable with finite mean \cite{B08}, this result is however difficult to apply in our setting as one need to find random variables that stochastically dominate the martingale difference sequence. We leave this question of extending the main result to unbounded $f$ as future work.
\end{rk}

As our first example to illustrate our main result Theorem \ref{thm:main}, we investigate the Jacobi process $X = \left(X_{t} \right)_{t\geq 0}$ on the state space $S=(0,1)$. The generator of this process is given by 
\begin{align}\label{eq:genJacobi}
\mathcal{A}^J=\left ( a-bx \right )\frac{d}{dx}+\frac{\sigma ^{2}}{2}x(1-x)\frac{d^2 }{d x^2},
\end{align}
where $a, b, \sigma \in \mathbb{R}$ are parameters of $X$ and are assumed to take on values such that $\alpha := \frac{2b }{\sigma ^{2}}-\frac{2a}{\sigma ^{2}}-1 > -1$ and $ \beta := \frac{2a }{\sigma ^{2}}-1 > -1$, i.e. $b > a > 0$ and $\sigma \in \mathbb{R}$. With these choices of parameters, the stationary distribution of Jacobi process is the Beta distribution with parameters $\alpha+1$ and $\beta+1$, where its density is governed by $\pi(x)=\frac{x^{\beta }(1-x)^{\alpha }}{B\left ( \alpha +1,\beta +1 \right )}$, and $B(\cdot,\cdot)$ denotes the Beta function. According to \cite[Appendix $B.3$]{AK07} and \cite[Section $2.1$]{FS08} $X$ is ergodic with $\sigma_{ess}(\mathcal{A}^J) = \emptyset$ and 
\begin{align*}
\lambda _{i}&= \frac{\sigma^{2} }{2}i\left ( i-1+\frac{2b}{\sigma ^{2}} \right ), \quad \sum_{i=1}^{\infty }\frac{1}{\lambda _{i}} = \frac{2}{\sigma ^{2}}\sum_{i=1}^{\infty }\frac{1}{i\left ( i-1+\frac{2b}{\sigma ^{2}} \right )}< \infty.
\end{align*} 
In view of Proposition \ref{prop:nscue} item \eqref{it:E}, $X$ is thus uniformly ergodic.
One major motivation for us to study such a process stems from its usage in financial modelling, where a more general form of Jacobi process has been employed to model exchange rates in a target zone, see \cite{LS07} and the references therein. In these models, one is often interested in the long-run average of the occupation time of the process in certain region $A$, say the occupation time of the exchange rate above or below a threshold. Unfortunately, distributional information on the functional $\int_0^t \1_{A}(X_s)\,ds$ is often inaccessible, where $\1_{A}$ is the indicator function of the set $A$. In practice, one may resort to the space average $\pi(\1_{A})$ as a natural approximation of the quantity of interest $\int_0^t \1_{A}(X_s)\,ds$, where the former is often easier to access than the latter. Our main result in Theorem \ref{thm:main} thus provides an invaluable tool and can be used to give non-asymptotic probabilistic error bounds on such approximation. Another situation where Theorem \ref{thm:main} is needed is about constructing confidence interval of the functional $\int_0^t \1_{A}(X_s)\,ds$. One can easily construct confidence band based on these large deviation probability. With these motivations in mind, we now apply Theorem \ref{thm:main} to the Jacobi process with $f = \1_A$ that gives:

\begin{corollary}\label{cor:main}
Suppose that $X = \left(X_{t} \right)_{t\geq 0}$ is the Jacobi process which is uniformly ergodic with generator given by \eqref{eq:genJacobi} and parameters $b > a > 0, \sigma \in \mathbb{R}$. For any $\varepsilon > 0$, $t > \frac{4t_{av}}{\varepsilon}$ and measurable subset $A \subseteq (0,1)$, we have
    \begin{align}\label{eq:jacobi}
    \mathbb{P}_{x}\left( \frac{1}{t}\int_{0}^{t}\1_{A}\left(X_{s}\right)ds-\pi\left(\mathbb \1_{A}\right)\geqslant \varepsilon  \right) \leqslant \exp\left \{ \frac{-2\left ( t\varepsilon -4t_{av} \right )^{2}}{(t+1)\left ( 4t_{av} +1\right )^{2}} \right \},
    \end{align}
    where
    $$t_{av} = \frac{2}{\sigma ^{2}}\sum_{i=1}^{\infty }\frac{1}{i\left ( i-1+\frac{2b}{\sigma ^{2}} \right )}$$
    is the average hitting time of the Jacobi process.
\end{corollary}

As our first remark, we note that the upper bound in \eqref{eq:jacobi} can be quite loose since it does not depend on the size of $A$. Such a bound indeed holds as long as we have $||f|| \leq 1$ in Theorem \ref{thm:main}. In addition, we see that this upper bound depends only on the parameters $b$ and $\sigma$ through $t_{av}$ but not on $a$. 

In our second example, we introduce the finite interval analogue of the Ornstein-Uhlenbeck process first studied by \cite{KS99} on the state space $S = (-\pi/2,\pi/2)$, where we take the drift to be $\mu(x) = -\rho \tan(x)$, the diffusion coefficient to be $\sigma^2(x) = 1$ and $\rho \geq 1/2$ to be a parameter. That is, the generator $\mathcal{A}^{O}$ is written as
\begin{align}\label{eq:AO}
\mathcal{A}^{O} = -\rho \tan(x) \dfrac{d}{dx} + \dfrac{1}{2} \dfrac{d^2}{dx^2}.
\end{align}
According to \cite[Section $2.1$]{FS08} $X$ is ergodic with $\sigma_{ess}(\mathcal{A}^O) = \emptyset$ and 
\begin{align}\label{eq:eigenOU}
\lambda _{i}&= i\left ( \rho + i/2 \right ), \quad \sum_{i=1}^{\infty }\frac{1}{\lambda _{i}} = \sum_{i=1}^{\infty }\frac{1}{i\left ( \rho + i/2 \right )}< \infty.
\end{align} 
In view of Proposition \ref{prop:nscue} item \eqref{it:E}, $X$ is thus uniformly ergodic for any $\rho \geq 1/2$. Specializing into the case $\rho = 1/2$, we see that the stationary distribution has density given by 
$$\pi(x) = \dfrac{\cos(x)}{2} \1_{x \in (-\pi/2,\pi/2)}.$$
For $u \in \mathbb{R}$, if we take $f(x) = e^{ux}\1_{x \in (-\pi/2,\pi/2)}$ with $||f|| \leq e^{u\pi/2}$ in Theorem \ref{thm:main}, the time integral becomes
$$\int_0^t f(X_s)\,ds = \int_{0}^{t} e^{uX_s}ds,$$
the exponential functional associated with $X$. Often distributional information of exponential functionals are difficult to obtain, see for instance the book \cite{Yor01}. One may approximate such functional by means of their space average $\pi(f)$, and our results come in handy since they give probabilistic error bound on such approximation. Theorem \ref{thm:main} now reads
\begin{corollary}\label{cor:main2}
	Suppose that $X = \left(X_{t} \right)_{t\geq 0}$ is the finite interval analogue of the Ornstein-Uhlenbeck process which is uniformly ergodic with generator given by \eqref{eq:AO} and parameter $\rho = 1/2$. For any $\varepsilon > 0$, $u \in \mathbb{R}$ and $t > \frac{4e^{u\pi/2}t_{av}}{\varepsilon}$, we have
	\begin{align}\label{eq:OU}
	\mathbb{P}_{x}\left( \int_{0}^{t} e^{uX_s}ds- \dfrac{2t\cosh(u\pi/2)}{1+u^2}   \geqslant t\varepsilon  \right) \leqslant \exp\left \{ \frac{-2\left ( t\varepsilon -4e^{u\pi/2}t_{av} \right )^{2}}{(t+1)e^{u\pi}\left ( 4t_{av} +1\right )^{2}} \right \},
	\end{align}
	where
	$$t_{av} = \sum_{i=1}^{\infty }\frac{2}{i\left ( i+1 \right )} = 2$$
	is the average hitting time of $X$.
\end{corollary}
For further concrete examples of uniformly ergodic diffusions with explicit eigenvalues information, we refer interested readers to the work of \cite{KS99,FS08}.

The rest of the paper is organized as follows. In Section \ref{sec:proofs}, we first present the proof of the main result Theorem \ref{thm:main}, followed by detailing the proof of Corollary \ref{cor:main} and the proof of Corollary \ref{cor:main2}.

\section{Proof of the main results}\label{sec:proofs}

\subsection{Proof of Theorem \ref{thm:main}}
Suppose without loss of generality that the mean of $f$ with respect to $\pi$ is zero, that is, $\pi(f)=0$.
To begin with, it follows readily from the induced operator norm of $Q^{\sharp}$ that we have
\begin{align}\label{eq:fhat}
\left \| \hat{f} \right \|\leqslant \left \| f \right \|\left \| Q^{\sharp} \right \|,
\end{align}
where we recall $\hat{f} = Q^{\sharp} f$ is the solution to the Poisson equation. Now, for the large deviation probability, we see that
\begin{align}  
          \mathbb{P}_{x}\left( \frac{1}{t}\int_{0}^{t}f(X_s)ds-\pi(f)\geqslant \varepsilon  \right)\label{eq: markov}
          &\leqslant e^{-\theta t\varepsilon }\mathbb{E}_{x}\left [ e^{\theta \int_{0}^{t}f\left ( X_{s} \right )ds} \right ]\\  
          &= e^{-\theta t\varepsilon }\mathbb{E}_{x}\left [ e^{-\theta \int_{0}^{t}\mathcal{A}\hat{f}\left ( X_{s} \right )ds} \right ].\label{eq: poi}
\end{align}
In the above equation, \eqref{eq: markov} comes from Markov inequality, which holds for any $\theta\geqslant 0$, while \eqref{eq: poi} follows from the Poisson equation $-\mathcal{A}\hat{f}=f$. Now, we explicitly construct a martingale that is useful in our analysis, namely
\begin{align}  
\mathcal{M}_{t}^{\hat{f}} := \hat{f}\left ( X_{t} \right )-\hat{f}\left ( X_{0} \right )-\int_{0}^{t}\mathcal{A}\hat{f}\left ( X_{s} \right )\,ds.\label{eq: mar}
\end{align}
Then by a classical result in \cite[Chapter $5$ Theorem $2.3$]{BW09}, we see that $\mathcal{M}_{t}^{\hat{f}}$ is a mean zero $\left \{ \mathcal{F}_{t}  \right \}$-martingale, where we again recall $ \mathcal{F}_{t} =\sigma \left \{ X_{u}:0\leqslant u\leqslant t \right \}$ is the filtration of $X$. Using \eqref{eq:fhat} and \eqref{eq: mar}, the tail bound in \eqref{eq: poi} is further upper bounded by
\begin{align}  
e^{-\theta t\varepsilon }\mathbb{E}_{x}\left [ e^{-\theta \int_{0}^{t}\mathcal{A}\hat{f}\left ( X_{s} \right )ds} \right ]
&= e^{-\theta t\varepsilon }\mathbb{E}_{x}\left [ e^{\theta\left ( \mathcal{M}_{t}^{\hat{f}}+\hat{f}\left ( X_{0} \right )-\hat{f}\left ( X_t \right )\right )    } \right ]   \nonumber \\
  &\leqslant e^{-\theta t\varepsilon }e^{2\theta \left \| Q^{\sharp} \right \|\left \| f \right \|} \mathbb{E}_{x}\left[e^{\theta \mathcal{M}_{t}^{\hat{f}} } \right].\label{eq:bound1}
\end{align}

Now, we proceed to examine the bound for $ \mathbb{E}_{x}\left[e^{\theta \mathcal{M}_{t}^{\hat{f}} } \right]$. In order to use the classical Hoeffding's lemma for bounded random variables \cite[Lemma $8.1$]{MR1}, we write $\mathcal{M}_{t}^{\hat{f}}$ as
\begin{align*} \mathcal{M}_{t}^{\hat{f}}&=\sum_{s=1}^{\left \lfloor t \right \rfloor}\left ( \mathcal{M}_{s}^{\hat{f}}-\mathcal{M}_{s-1}^{\hat{f}} \right )+\mathcal{M}_{t}^{\hat{f}}-\mathcal{M}_{\left \lfloor t \right \rfloor}^{\hat{f}}.
\end{align*}
As a result, to bound the martingale $\mathcal{M}^{\hat{f}}_t$ it suffices to bound the martingale differences $\mathcal{M}^{\hat{f}}_s - \mathcal{M}^{\hat{f}}_{s-1}$. Using the definition of $\mathcal{M}_t^{\hat{f}}$ in \eqref{eq: mar}, these bounds are given by, for $s=1,2,...\left \lfloor t \right \rfloor$, 
 \begin{align}
 \mathcal{M}_{s}^{\hat{f}}-\mathcal{M}_{s-1}^{\hat{f}}&=\hat{f}\left ( X_{s} \right )-\hat{f}\left ( X_{s-1} \right )-\int_{s-1}^{s}\mathcal{A}\hat{f}\left ( X_{s} \right )\,ds \nonumber\\
 \left |\mathcal{M}_{s}^{\hat{f}}-\mathcal{M}_{s-1}^{\hat{f}} \right| \nonumber
&\leqslant 2\left \| \hat{f} \right \|+\left| \int_{s-1}^{s}f\left ( X_{s} \right )\,ds \right| \\\nonumber
&\leqslant 2\left \|  Q^{\sharp} \right \|\left \| f \right \|+\left \| f \right \|\\
&=\left ( 2\left \|  Q^{\sharp} \right \| +1\right )\left \| f \right \| \label{eq:part1},
\end{align} 
where we use the Poisson equation in the first inequality and \eqref{eq:fhat} in the second inequality. Similarly, 
\begin{align}
\mathcal{M}_{t}^{\hat{f}}-\mathcal{M}_{\left \lfloor t \right \rfloor}^{\hat{f}}\nonumber
&=\hat{f}\left ( X_{t} \right )-\hat{f}\left ( X_{\left \lfloor t \right \rfloor} \right )-\int_{\left \lfloor t \right \rfloor}^{t}\mathcal{A}\hat{f}\left ( X_{s} \right )\,ds\\
\left|\mathcal{M}_{t}^{\hat{f}}-\mathcal{M}_{\left \lfloor t \right \rfloor}^{\hat{f}}  \right|
 &\leqslant 2\left \| \hat{f} \right \|+\left \| f \right \| \leq \left ( 2\left \|  Q^{\sharp} \right \|+1 \right )\left \| f \right \|.\label{eq:part2}
\end{align}
It follows from double expectation, \eqref{eq:bound1}, \eqref{eq:part1} and \eqref{eq:part2} that the upper bound in \eqref{eq: poi} becomes
\begin{align}  
          \mathbb{E}_{x} \left[ e^{ \theta\left ( \sum_{s=1}^{\left \lfloor t \right \rfloor}\left ( \mathcal{M}_{s}^{\hat{f}}-\mathcal{M}_{s-1}^{\hat{f}} \right )+\mathcal{M}_{t}^{\hat{f}}-\mathcal{M}_{\left \lfloor t \right \rfloor}^{\hat{f}}\right )} \right]  
          &=\mathbb{E}_{x}\ \left [  e ^{ \sum_{s=1}^{\left \lfloor t \right \rfloor}\left ( \mathcal{M}_{s}^{\hat{f}}-\mathcal{M}_{s-1}^{\hat{f}} \right )}\mathbb{E}_{x}\ \left [  e ^{ \theta \left ( \mathcal{M}_{t}^{\hat{f}}-\mathcal{M}_{\left \lfloor t \right \rfloor}^{\hat{f}} \right )}\bigg| \mathcal{F}_{\left \lfloor t \right \rfloor}\right]\right]   \nonumber\\
&\leqslant \mathbb{E}_{x}\ \left [  e ^{ \theta \sum_{s=1}^{\left \lfloor t \right \rfloor}\left ( \mathcal{M}_{s}^{\hat{f}}-\mathcal{M}_{s-1}^{\hat{f}} \right )}\right ] e^{\frac{\theta ^{2}\left ( 2\left \| Q^{\sharp} \right \|+1 \right )^{2}\left \| f \right \|^{2}}{8}}\nonumber   \\
&\leqslant e^{\frac{\theta ^{2}(t+1)\left ( 2\left \| Q^{\sharp} \right \|+1 \right )^{2}\left \| f \right \|^{2}}{8}}\nonumber, 
\end{align}
where the first and second inequality follows from repeated applications of the Hoeffding's lemma \cite[Lemma $8.1$]{MR1}. Finally, collecting the above results the tail bound is given by
\begin{align}
\mathbb{P}_{x}\left( \frac{1}{t}\int_{0}^{t}f(X_s)ds-\pi(f) \geqslant \varepsilon  \right)
&\leqslant \exp\left ( \frac{\theta ^{2}(t+1)\left ( 2\left \| Q^{\sharp} \right \|+1 \right )^{2}\left \| f \right \|^{2}}{8} -\theta t\varepsilon +2\theta \left \| Q^{\sharp} \right \|\left \| f \right \|\right ),\label{eq:f0}
\end{align}
which is minimized at $\theta = \theta^*$ where 
\begin{align}
\theta^*&=\frac{4t\varepsilon -8\left \| Q^{\sharp} \right \|\left \| f \right \|}{(t+1)\left ( 2\left \| Q^{\sharp} \right \|+1 \right )^{2}\left \| f \right \|^{2}}\nonumber.
\end{align}
Desired result follows by substituting $\theta = \theta^*$ into \eqref{eq:f0}.

\subsection{Proof of Corollary \ref{cor:main}}
Desired result follows from taking $f= \1_{A}$ in Theorem \ref{thm:main} and utilizing the follow bound on the induced operator norm of the deviation kernel $Q^{\sharp}$:
\begin{align}
\left \| Q^{\sharp} \right \|\leqslant 2t_{av},
\end{align} 
see e.g. \cite[Theorem $1.1$]{Choi18}. As for the expression of the average hitting time $t_{av}$, the eigentime identity \cite{CM15} gives 
\begin{align}
t_{av}=\sum_{i=1}^{\infty }\frac{1}{\lambda _{i}} = \frac{2}{\sigma ^{2}}\sum_{i=1}^{\infty }\frac{1}{i\left ( i-1+\frac{2b}{\sigma ^{2}} \right )}< \infty,
\end{align}
where $0 = \lambda_0  < \lambda_1 \leq \lambda_2 \leq \ldots$ are the eigenvalues of $-\mathcal{A}^J$ which are given by, for $i = 0,1,2,\ldots$,
\begin{align*}
\lambda _{i}&= \frac{\sigma^{2} }{2}i\left ( i-1+\frac{2b}{\sigma ^{2}} \right ),
\end{align*}
see e.g. \cite[Appendix $B.3$]{AK07}.

\subsection{Proof of Corollary \ref{cor:main2}}
Desired result follows from taking $f(x) = e^{ux}\1_{x \in (-\pi/2,\pi/2)}$ in Theorem \ref{thm:main}, and using
$$\pi(f) = \dfrac{1}{2} \int_{-\pi/2}^{\pi/2} e^{ux} \cos(x) \,dx = \dfrac{2 \cosh(u\pi/2)}{1+u^2}$$
as well as the following bound on the induced operator norm of the deviation kernel $Q^{\sharp}$:
\begin{align*}
\left \| Q^{\sharp} \right \|\leqslant 2t_{av} = \sum_{i=1}^{\infty }\frac{4}{i\left ( i+1 \right )} = 4,
\end{align*} 
where again the first equality follows from \cite{CM15} with $\lambda_i$ being given in \eqref{eq:eigenOU} with parameter $\rho = 1/2$.
\section*{Acknowledgement}
\noindent \textbf{Acknowledgement}.
We thank the anonymous referee for constructive comments that improve the presentation of the manuscript. This work is partially supported by the Chinese University of Hong Kong, Shenzhen grant PF01001143.

\bibliographystyle{abbrvnat}
\bibliography{library}
\end{document}